\documentclass{article}
\usepackage{graphicx} 
\usepackage{multirow}
\usepackage{latexsym} 
\usepackage{amsmath} 
\usepackage{amsfonts} 
\usepackage{amssymb} 
\usepackage{soul} 
\usepackage{ulem} 
\usepackage{lscape} 
\begin{document}     
\baselineskip 7.8mm
\baselineskip 5.5mm

\title{\bf 
Fourier and Zak transforms of multiplicative characters \\
}  

\vspace{4mm}

\author{ 
\\ Andrzej K. Brodzik$^{1,2}$  and Richard Tolimieri$^{3}$ \\   \\  
{\small $^{1}$Boston University, 
\small Department of Psychological and Brain Sciences, 
Boston, MA 02215, USA} \\
\small $^{2}$Gda$\rm\acute{n}$sk University of Technology, Faculty of Electronics, Telecommunications and Informatics,
Gda$\acute{\rm n}$sk, Poland \\
\small email: andrzej.k.brodzik@gmail.com \\ \\
\small $^{3}$Prometheus Inc.,
Newport, RI 02840, USA \\
\small email: richard@psypher.org \\ \\
}

\date{}

\maketitle


\baselineskip 5.0mm
\baselineskip 4.9234050mm

\begin{center}
\section*{\small Abstract}
\end{center}
In this paper we derive formulas for the $N$-point discrete Fourier transform and the $R_1\times R_2$ finite Zak transform 
of multiplicative characters on $\Bbb{Z}/N$, 
where $N$ is an odd integer,
and $R_1$ and $R_2$ are co-prime factors of $N$.
In one special case this permits computation of the discrete Fourier transform and the finite Zak transform 
of the Jacobi symbol, the modified Jacobi sequence, and the Golomb sequence.
In other cases, not addressed here, 
this permits computation of the discrete Fourier transform and the finite Zak transform 
of certain complex-valued sequences.
These results 
constitute, to our knowledge, 
the first unified treatment of key Fourier and Zak space properties of multiplicative characters.
These results also provide a convenient framework for 
the design of new character-based sequences.
\\ 

\noindent
{\bf Index terms:} 
Chinese Remainder Theorem, 
Bj$\rm\ddot{o}$rck sequences, 
discrete Fourier transform,
finite Zak transform, 
Gauss sum, Golomb sequences, Good-Thomas Prime Factor algorithm,
ideal sequences, Jacobi symbol, Legendre symbol, modified Jacobi sequences, 
multiplicative characters, Quadratic Reciprocity Law, 
twin-prime sequences. \\
\\

\noindent
\section{Introduction}

\noindent
The focus of this work is on Fourier properties of multiplicative characters.
A multiplicative character on $\Bbb{Z}/N$ is a group homomorphism from its group of units, $U(N)$, to the multiplicative group of
non-zero complex numbers, $\Bbb{C}^{\ast}$ [20].
Multiplicative characters are of theoretical interest, as they link concepts in group theory, linear algebra, harmonic analysis, 
and number theory. 
Among multiplicative characters a class of multiplicative characters, 
the primitive multiplicative characters, 
is especially important, for two main reasons.
First, the discrete Fourier transform (DFT) of
primitive multiplicative characters
has a simple structure explicitly involving the Gauss sums 
that can be used to re-derive certain classical number-theoretical results.
Second, and of more direct interest in this work, primitive multiplicative characters 
are closely associated with sequences with favorable pseudo-random and correlation properties, including
difference sets and quadratic residue sequences. 
These sequences are often used in cryptography [25], [28],
communications [16], [23] and radar [16], [23], [27-28]. 
The association between characters and sequences has been explored previously, 
among others by Turyn in [32] (in the construction of Hadamard difference sets) 
and by Scholtz and Welch in [27] (in investigation of sequences with good correlation properties). \\


\noindent
The principal results of this work include: 

\begin{enumerate}
\item
The formulas for the $N$-point DFT and the $R_1\times R_2$ 
finite Zak transform (FZT) 
of primitive multiplicative characters on $\Bbb{Z}/N$,
where $N=R_1R_2$ is an arbitrary odd integer and $(R_1,R_2)=1$.
While the DFT formula was given previously in [4], 
the presentation given here is more pedagogical, as it builds on detailed analysis
of special cases, $\Bbb{Z}/p$ and $\Bbb{Z}/p^m$, and as it provides a link to a subsequent time-frequency analysis
of sequences.
\item
The use of the formulas for the DFT and the FZT of primitive characters leads to the derivation of the DFT and the FZT
of three special sequences: the Jacobi symbol, the modified Jacobi sequence and the Golomb sequence.
The FZT result appears especially interesting, as it suggests a link 
between the constructions of complex-valued and binary sequences. 
\item
The use of an intermediate result in the derivation of the DFT of the Jacobi symbol that links the Jacobi symbol
with Gauss sums (Theorem 6) and a prior result for the trace of the DFT matrix leads to 
re-derivation of the Quadratic Reciprocity Law (QRL).
\end{enumerate}

\noindent
The DFT of primitive multiplicative characers is derived twice, in two different contexts and using different tools. 
The first derivation relies, in part, on prior results, summarized in Section 2, for 
the
primitive multiplicative characters
on $\Bbb{Z}/p$ 
and for the primitive multiplicative characters
on $\Bbb{Z}/p^m, \ m\ge 2$, 
and, in part, on the Chinese Remainder Theorem (CRT).
The CRT
defines a ring isomorphism between 
$\Bbb{Z}/N$ and $\Bbb{Z}/R_1\times \Bbb{Z}/R_2$. 
This ring isomorphism establishes a relationship between
the multiplicative characters on   
$\Bbb{Z}/N$, and those on $\Bbb{Z}/R_1$ and $\Bbb{Z}/R_2$,
\begin{eqnarray*}
\chi(a_1e_1+a_2e_2)=\chi_1(a_1)\chi_2(a_2), \ \ a_1\in U(R_1), \ \ a_2\in U(R_2),
\end{eqnarray*}
where $e_1$ and $e_2$ are idempotents for $\Bbb{Z}/N$ corresponding to the factorization $N=R_1R_2$.
This relationship permits 
construction of the
primitive multiplicative characters
on $\Bbb{Z}/N$ from 
the primitive multiplicative characters
on $\Bbb{Z}/R_1$ and $\Bbb{Z}/R_2$ (Theorem 3), and
computation of the DFT
of primitive multiplicative characters on $\Bbb{Z}/N$ (Theorem 4). The last result is given by the formula
\begin{eqnarray*}
\hat{\chi}(a)=G(\chi)\chi^{\ast}(a), 
\end{eqnarray*}
where
\begin{eqnarray*}
G(\chi)=G(\chi_1)G(\chi_2)\chi_1^{\ast}(e_1/R_2)\chi_2^{\ast}(e_2/R_1). 
\end{eqnarray*}
In Section 4 we present a second proof of Theorem 4, based on the Good-Thomas Prime Factor (GTPF) algorithm. 
This derivation is 
more explicit, as it relies on the relationship between the DFT matrix $F(N)$ and the tensor product of the 
associated DFT matrices $F(R_1)$ and $F(R_2)$, perturbed by permutation matrices. \\ 

\noindent
The second main theoretical result of this work (Theorem 9) is the $R_1\times R_2$ FZT of a primitive multiplicative character on $\Bbb{Z}/N$,
\begin{eqnarray*}
X_{R_2}(j,k)=G(\chi_2)
\chi_2(R_1)e^{-2\pi i \frac{R_1^{-1}jk}{R_2}}
\chi_2^{\ast}(j)
\chi_1(k), \ \ \
0\leq j<R_2, \ \ 0\leq k<R_1. 
\end{eqnarray*} 
This result can be viewed as a generalization of the 
result for the DFT. 
The formula for the FZT is more broadly applicable than the formula for the DFT, in the sense that it does not require 
both multiplicative characters on $\Bbb{Z}/R_1$ and $\Bbb{Z}/R_2$ to
be primitive.  
The use of the FZT in the analysis of multiplicative characters is new; it complements and sometimes is preferable 
to the use of the DFT, especially in sequence design applications. 
We elaborate on this in Section 7.  \\ 

\noindent


\noindent
One consequence of these results is that they avail convenient tools
for investigating special cases that are important in sequence design applications.
We focus here, in particular, on three 
sequences:
the Jacobi symbol, 
the modified Jacobi sequence [22], and the Golomb sequence [15].
The Jacobi symbol is the second order primitive multiplicative character and a well-known ternary sequence.
The last two sequences 
are not multiplicative characters, but they can be obtained from multiplicative characters by simple modifications. 
The importance of the modified Jacobi sequence stems from the fact that while, 
like the Jacobi symbol, it has good aperiodic correlation properties and large linear complexity,
unlike the Jacobi symbol, it is binary.
The importance of the Golomb sequence stems from the fact that,
unlike the Jacobi symbol and the modified Jacobi sequence, it is a 
constant magnitude ideal sequence, that is a sequence with an ideal autocorrelation. \\ 

\noindent
We derive the DFT (Corollary 4 and Theorems 7 and 8) and the FZT (Corollary 5 and Theorems 10 and 11) 
of all three sequences. 
The DFT results are not entirely new; partial results for all three cases have appeared before.
The DFT of the Jacobi sequence was used by Jensen {\it et al} in [22], to facilitate computation of merit factor of binary sequences, 
but the formula is not stated there explicitly in a complete, general form.
The DFT of the modified Jacobi sequence was given in [22] as well, but the formula is more complex than
necessary, and it depends on certain 
auxiliary parameters.
The DFT of the Golomb sequence for the special case of prime length (referred there as the Bj$\rm\ddot{o}$rck sequence) 
was given by Popovic in [24].
Our results generalize and unify these results.
In further research the general result for the DFT of multiplicative characters 
can be used to identify new sequences and to compute their DFTs.
The FZT results are entirely new.
Among others, these results show that the FZT of the Golomb sequence, 
while, unlike the FZT of a chirp [10], is neither unimodular nor semi-unimodular,
it is highly structured, suggesting Zak space sequence design approaches, similar to those developed for chirps [2], [6-11], 
might be viable. \\

\noindent

\noindent
This work is a sequel to a sequence of papers and books.
In particular, it is a sequel to [10], which, in part, synthesizes our principal results for chirps.
The DFT formula in Theorem 4 extends the result in [31] for the DFT of multiplicative characters 
on $\Bbb{Z}/p^m$ to the general case, and is
an analogue of the formulas for the DFT of a discrete chirp given by Theorem 1 in [6] and Theorem 1 in [10].
Associated with the DFT formula in Theorem 4, the FZT formula in Theorem 7  
is an analogue of the formulas for the FZT of a discrete chirp given by Theorem 5 in [7] and Theorems 4 and 5 in [10]. \\

\noindent
The content of this paper is as follows.
In Section 2 we review prior results on the DFT of multiplicative characters on $\Bbb{Z}/N$ for $N=p^m$, $m\ge 1$.
In Section 3 we use the CRT to extend the results of Section 2 to the case of 
$N=R_1R_2$, $N$ odd, $(R_1,R_2)=1$.
In Section 4 we re-derive this result using the GTPF algorithm.
In Section 5 we derive a formula for the DFT of the Jacobi symbol. 
In Section 6 we derive formulas for the DFT of the modified Jacobi sequence and the Golomb sequence.
In Section 7 we derive formulas for the FZT of multiplicative characters, 
the Jacobi symbol, the modified Jacobi sequence, and the Golomb sequence.
In Section 8 we discuss further work.
In Appendix 1 we summarize basic properties of the Legendre and Jacobi symbols.
In Appendix 2 we re-derive the QRL using Theorem 6 and a prior result for the trace of the DFT matrix.

\vspace{2mm}

\newpage

\noindent
\section{DFT of multiplicative characters in a special case}


\noindent
We will start by introducing basic notation.
$N$ is an odd positive integer with distinct prime power factorization 
$N=p_1^{c_1}...p_r^{c_r}$, 
${\Bbb{Z}}/N$
is the ring of integers modulo $N$ and $U(N)$ is its multiplicative group of units.
$U(N)$ consists of all $a\in {\Bbb{Z}}/N$ such that $(a,N)=1$ and has order
\begin{eqnarray}
t=p_1^{c_{1}-1}...p_r^{c_{r}-1}(p_1-1)...(p_r-1). 
\end{eqnarray}
It follows that the order of $U(N)$ is even.   \\

\noindent
Suppose $x$ is an arbitrary $N$-periodic sequence. 
The $N$-point Fourier transform of $x$, $\hat{x}$, is given by 
\begin{eqnarray}
\hat{x}(a)=\sum_{b\in \Bbb{Z}/N} x(b)e^{2\pi i\frac{ab}{N}}, \ \ a\in {\Bbb{Z}}/N.
\end{eqnarray} 

\vspace{3mm}

\noindent
A multiplicative character on $\Bbb{Z}/N$ 
is a group homomorphism, 
$\chi: U(N)\longrightarrow {\Bbb{C}}^{\ast}$, 
where ${\Bbb{C}}^{\ast}$
is the multiplicative group of non-zero complex numbers, by definition, 
\begin{eqnarray}
\chi(ab)=\chi(a)\chi(b), \ \ a, b\in U(N).
\end{eqnarray}
Multiplicative characters take their values on the $t$-th root of unity, $t$ the order of $U(N)$.
The set of multiplicative characters forms a group of order $t$ under the rule
\begin{eqnarray}
\chi_1\chi_2(a)=\chi_1(a)\chi_2(a), \ \ a\in U(N), \ \ \chi_1, \chi_2 \ \ {\rm multiplicative \ characters \ on} \ \Bbb{Z}/N.
\end{eqnarray}
In particular
\begin{eqnarray}
\chi(a^{-1})=\chi^{\ast}(a)=\chi^{-1}(a), \ \ a\in U(N). 
\end{eqnarray}
We set the values of multiplicative characters equal to zero off of $U(N)$ and assume that they are defined on $\Bbb{Z}$ by
periodicity mod $N$,    
\begin{eqnarray}
\chi(a)=\chi(a_1), \ \ a\in \Bbb{Z}, \ \ a_1\equiv a \!\!\mod N. 
\end{eqnarray}

\noindent
An important role in this work is played by the properties of $U(p)$ and its multiplicative characters, $p$ an odd prime.
$U(p)$ is a cyclic group of order $t=p-1$.
Denote a generator of $U(p)$ by $g$, so that $U(p)$ consists of the powers 
\begin{eqnarray}
g^a, \ \ 0\leq a<t.
\end{eqnarray} 
There are $t$ distinct multiplicative characters on $\Bbb{Z}/p$ defined by
\begin{eqnarray}
\chi_j(g^a)=e^{2\pi i\frac{ja}{t}}, \ \ 0\leq j<t, \ \ 0\leq a<t.
\end{eqnarray} 

\noindent
By well-known results (see, for example, [20])
\begin{eqnarray}
\hat{\chi}_0(a)=
\begin{cases}
\ p-1,     
&  a=0, \\   
\ -1,     
&  otherwise, \\   
\end{cases}
\end{eqnarray} 
and
\begin{eqnarray}
\hat{\chi}_j(a)=\hat{\chi}_j(1)\chi_j^{\ast}(a), \ \ \ j>0, \ \ \ a\in \Bbb{Z}/p, 
\end{eqnarray} 
where
\begin{eqnarray}
\hat{\chi}_j(1)\neq 0. 
\end{eqnarray} 
In fact, for $j>0$,
\begin{eqnarray}
|\hat{\chi}_j(1)|^2=p. 
\end{eqnarray} 


\vspace{3mm}


\noindent
Suppose $\chi$ is an arbitrary multiplicative character on $U(N)$ and $\hat{\chi}$ is its $N$-point Fourier transform
\begin{eqnarray}
\hat{\chi}(a)=\sum_{b\in U(N)} \chi(b)e^{2\pi i\frac{ab}{N}}, \ \ a\in {\Bbb{Z}}/N.
\end{eqnarray}
For $a\in U(N)$
\begin{eqnarray}
\chi(a)\hat{\chi}(a)=
\sum_{b\in U(N)} \chi(ab)e^{2\pi i\frac{ab}{N}} 
=\sum_{c\in U(N)} \chi(c)e^{2\pi i\frac{c}{N}}=\hat{\chi}(1). 
\end{eqnarray}

\noindent
This leads to the following well-known result [20].  \\

\noindent
{\bf Theorem 1} \ {\it 
If $\chi$ is a multiplicative character on $\Bbb{Z}/N$, then} 
\begin{eqnarray}
\hat{\chi}(a)=\hat{\chi}(1)\chi^{\ast}(a), \ \ a\in U(N).
\end{eqnarray}

\vspace{2mm}

\noindent
There are two mutually exclusive cases. If
\begin{eqnarray}
\hat{\chi}(1)=0, 
\end{eqnarray}
then $\hat{\chi}$ vanishes on $U(N)$. Otherwise $\hat{\chi}$ is a non-zero scalar multiple of $\chi^{\ast}$ on $U(N)$. 
In this case we set
\begin{eqnarray}
G(\chi)=\hat{\chi}(1), 
\end{eqnarray}
and call $G(\chi)$ the Gauss sum of the multiplicative character $\chi$. \\

\noindent
We call a multiplicative character $\chi$ on $\Bbb{Z}/N$ {\it primitive} if
\begin{eqnarray}
\hat{\chi}(a)=G(\chi)\chi^{\ast}(a), \ \ a\in {\Bbb{Z}}/N.
\end{eqnarray}
By necessity, if $\chi$ is primitive, then $\hat{\chi}(1)\neq 0$, since otherwise $\hat{\chi}$ vanishes everywhere,
and therefore $\chi$ must vanish everywhere. \\

\noindent
If $N=p$, the multiplicative characters $\chi_j,  \ 0<j<p-1$, are primitive, while $\chi_0$ is not,
since $\hat{\chi}_0(0)\neq 0$.
In the next section we construct the set of all primitive multiplicative characters on $\Bbb{Z}/N$
for arbitrary $N$. 
We remark that Theorem 1 
does not address the values of $\hat{\chi}(a),  \ a\notin 
U(N)$, 
a problem that has been studied in [31]. \\

\noindent
Next, we will describe the
primitive multiplicative characters on $\Bbb{Z}/p^m, \  m\geq 2$. 
Suppose $\chi$ is a multiplicative character on $\Bbb{Z}/p^m,  \ m\geq 2$,
satisfying $\hat{\chi}(1)\neq 0$. 
By Theorem 1,
\begin{eqnarray}
\hat{\chi}(a)=G(\chi)\chi^{\ast}(a), \ \ a\in U(p^m).
\end{eqnarray}
Below we will show that $\hat{\chi}$ vanishes outside 
$U(p^m)$.
Take $a\notin U(p^m)$ and write
\begin{eqnarray}
a=pa', \ \ 0\leq a'<p^{m-1}.
\end{eqnarray}
Consider the subgroup $\Delta$ of $U(p^m)$,
\begin{eqnarray}
\Delta=1+p^{m-1}{\Bbb{Z}}/p^m.
\end{eqnarray}
If $\chi(\Delta)=1$,
then for $v\in U(p^m)$ 
\begin{eqnarray}
\chi(v+p^{m-1}{\Bbb{Z}})=\chi(v)
\chi(1+v^{-1}p^{m-1}{\Bbb{Z}})=\chi(v)
\end{eqnarray}
and
$\chi$ is periodic mod 
$p^{m-1}{\Bbb{Z}}$. It follows that
 $\hat{\chi}$ is decimated mod $p\Bbb{Z}$, contradicting the assumption $\hat{\chi}(1)\neq 0$. 
It follows that we can find $c\in\Delta$, such that 
$\chi(c)\neq 1$. 
Since $pc\equiv p \!\! \mod p^m$, we have 
\begin{eqnarray}
w^{ab}=w^{pa'b}=w^{pca'b}=w^{abc}, \ \ w=e^{2\pi i\frac{1}{N}}, \ \ b\in U(p^m),
\end{eqnarray}
and
\begin{eqnarray}
\hat{\chi}(a)=\sum_{b\in U(p^m)}\chi(b)w^{abc}=\chi(c)\hat{\chi}(a). 
\end{eqnarray}
Since
$\chi(c)\neq 1$, $\hat{\chi}(a)=0$, proving (18) holds for $a\notin U(p^m)$. This yields the next result [31]. \\

\noindent
{\bf Theorem 2} \ {\it 
If $\chi$ is a  multiplicative character on $\Bbb{Z}/p^m, \ \ m\geq 2$,
satisfying the condition $\hat{\chi}(1)\neq 0$, then $\chi$ is primitive.} \\

\noindent
We can explicitly describe the primitive multiplicative characters
on $\Bbb{Z}/p^m, \ m\geq 2$ [31].
For an odd prime $p$, $U(p^m)$ is a cyclic group of order $t=p^{m-1}(p-1)$.
Suppose $z$ is a generator of $U(p^m)$.
The multiplicative characters on
$\Bbb{Z}/p^m$ are completely determined by their values on generator $z$ and are given by
\begin{eqnarray}
\chi_l(z)=e^{2\pi i\frac{l}{t}}, \ \ \ 0\leq l<t.
\end{eqnarray}
The primitive multiplicative characters are given by
\begin{eqnarray}
\chi_l, \ \ \ 0< l<t, \ \ \ (l, p^m)=1.
\end{eqnarray}
The non-primitive multiplicative characters are $\chi_0$ and
\begin{eqnarray}
\chi_l, \ \ \ 0< l<t, \ \ \ (l, p^m)>1.
\end{eqnarray}

\vspace{3mm}


\section{DFT of multiplicative characters and the CRT}

In this section we use the CRT to extend the results of Section 2 
to an arbitrary integer $N=R_1R_2$, $(R_1,R_2)=1$. \\

\noindent
Suppose $N=R_1R_2$. 
The mapping
\begin{eqnarray}
\Phi (a)= (a_1, a_2),  \ \ \ a_1\equiv a \!\!\! \mod R_1, \ \ a_2\equiv a \!\!\! \mod R_2,  
\end{eqnarray}
is a ring homomorphism of 
${\Bbb{Z}}/N$ into 
${\Bbb{Z}}/R_1\times {\Bbb{Z}}/R_2$.
The CRT states that if $(R_1, R_2)=1$, then $\Phi$ is a ring isomorphism of 
${\Bbb{Z}}/N$ onto 
${\Bbb{Z}}/R_1\times {\Bbb{Z}}/R_2$. \\

\noindent
Throughout $(R_1, R_2)=1$. 
By one form of the CRT there exist
\begin{eqnarray}
e_1, e_2 \in {\Bbb{Z}}/N 
\end{eqnarray}
such that
\begin{eqnarray}
e_1\equiv 1 \!\!\! \mod R_1 , \ \ \ 
e_1
\equiv 0 \!\!\! \mod R_2 \ \ 
\end{eqnarray}
and
\begin{eqnarray}
e_2\equiv 0 \!\!\! \mod R_1 , \ \ \ 
e_2
\equiv 1 \!\!\! \mod R_2. \ \ 
\end{eqnarray}
$e_1$ and $e_2$ satisfy the idempotent conditions
\begin{eqnarray}
e_1^2\equiv e_1 \!\!\! \mod N, \ \ \ 
e_2^2\equiv e_2 \!\!\! \mod N, \ \ \ 
e_1e_2\equiv 0 \!\!\! \mod N,  \ \ \ 
e_1+e_2\equiv 1 \!\!\!\mod N. 
\end{eqnarray}


\noindent
The mapping
\begin{eqnarray}
\bar{\Phi} (a_1, a_2)\equiv a_1e_1+a_2e_2 \!\!\! \mod N, \ \ \ a_1\in\Bbb{Z}/R_1, \ \ a_2\in\Bbb{Z}/R_2  
\end{eqnarray}
is a ring isomorphism from $\Bbb{Z}/R_1\times \Bbb{Z}/R_2$ onto $\Bbb{Z}/N$, providing the inverse to $\Phi$. 
In particular, every $a$ in $\Bbb{Z}/N$ can be uniquely written as
\begin{eqnarray}
a= a_1e_1+a_2e_2 \!\!\! \mod N,  \ \ \ a_1\in \Bbb{Z}/R_1, \ \ a_2\in \Bbb{Z}/R_2. 
\end{eqnarray} 

\noindent
$U(N)$ consists of the points
\begin{eqnarray}
 a_1e_1+a_2e_2 \!\!\! \mod N,  \ \ \ a_1\in U(R_1), \ \ a_2\in U(R_2). 
\end{eqnarray} 

\vspace{2mm}

\noindent
Suppose $\chi_1$ and $\chi_2$ are multiplicative characters 
on $\Bbb{Z}/R_1$ and on $\Bbb{Z}/R_2$, respectively, then
\begin{eqnarray}
 \chi(a_1e_1+a_2e_2)=\chi_1(a_1)\chi_2(a_2), \ \ \ a_1\in U(R_1), \ \ a_2\in U(R_2), 
\end{eqnarray}
is a multiplicative character on $\Bbb{Z}/N$. In the inverse direction, if $\chi$ is a multiplicative character on $\Bbb{Z}/N$, then
\begin{eqnarray}
\chi_1(a_1)= \chi(a_1e_1+e_2) \ \ \ {\rm and} \ \ \
\chi_2(a_2)= \chi(e_1+e_2a_2),
\end{eqnarray}
are multiplicative characters
on $\Bbb{Z}/R_1$ and on $\Bbb{Z}/R_2$. \\

\noindent
Also, 
\begin{eqnarray}
\chi_1(a_1) 
\chi_2(a_2)= 
\chi(a_1e_1+e_2)
\chi(e_1+e_2a_2)=
\chi(a_1e_1+e_2a_2).
\end{eqnarray} 
In this way we can identify multiplicative characters on $\Bbb{Z}/N$ with pairs of multiplicative characters
on $\Bbb{Z}/R_1$ and on $\Bbb{Z}/R_2$. 
We can also write
\begin{eqnarray}
\chi(a)=\chi_1(a_1) 
\chi_2(a_2), 
\ \ a_1\equiv a \!\!\! \mod R_1, 
\ \ a_2\equiv a \!\!\! \mod R_2. 
\end{eqnarray}  
Setting
\begin{eqnarray}
a\equiv a_1e_1+a_2e_2 \!\!\! \mod N 
\ \ \ {\rm and} \ \ \
b\equiv b_1e_1+b_2e_2 \!\!\! \mod N. 
\end{eqnarray} 
in
\begin{eqnarray}
\hat{\chi}(a)=\sum_{b\in U(N)} \chi(b)e^{2\pi i\frac{ab}{N}}, 
\end{eqnarray}
we obtain
\begin{eqnarray}
\hat{\chi}(a_1e_1+a_2e_2)=
\sum_{b_1\in U(R_1)} \chi_1(b_1)e^{2\pi i\frac{e_1}{R_2}\frac{a_1b_1}{R_1}} 
\sum_{b_2\in U(R_2)} \chi_2(b_2)e^{2\pi i\frac{e_2}{R_1}\frac{a_2b_2}{R_2}} 
= \hat{\chi}_1\left(a_1\frac{e_1}{R_2}\right) 
 \hat{\chi}_2\left(a_2\frac{e_2}{R_1}\right). 
\end{eqnarray}

\noindent
Since $e_1\equiv 0\mod R_2$ and $e_1\equiv 1\mod R_1$,
we have that $e_1/R_2$ is an integer relatively prime to $R_1$.
In the same way $e_2/R_1$ is an integer relatively prime to $R_2$, leading to the next result. \\

\noindent
{\bf Theorem 3} \ {\it 
If $\chi_1$ and $\chi_2$ are multiplicative characters on $\Bbb{Z}/R_1$ and on $\Bbb{Z}/R_2$, and
$\chi$ is the corresponding multiplicative character on $\Bbb{Z}/N$, then
\begin{eqnarray}
\hat{\chi}(a_1e_1+a_2e_2)=\hat{\chi}_1\left(a_1\frac{e_1}{R_2}\right) 
\hat{\chi}_2\left(a_2\frac{e_2}{R_1}\right), 
\ \ \ a_1\in {\Bbb{Z}}/R_1, 
\ \ a_2\in {\Bbb{Z}}/R_2, 
\end{eqnarray} 
where
$\frac{e_1}{R_2}\in U(R_1)$ and $\frac{e_2}{R_1}\in U(R_2)$. } \\

\noindent
Since $e_1+e_2\equiv 1\mod N$,
\begin{eqnarray}
\hat{\chi}(1)=
\hat{\chi}_1\left(\frac{e_1}{R_2}\right)
\hat{\chi}_2\left(\frac{e_2}{R_1}\right),
\end{eqnarray}
proving the following. \\

\noindent
{\bf Corollary 1} \ {\it
$\hat{\chi}(1)\neq 0$ if and only if 
$\hat{\chi}_1(1)\neq 0$ and 
$\hat{\chi}_2(1)\neq 0$. }  \\ 

\noindent
Moreover, since 
$\hat{\chi}(a)= 0$ for 
$a\notin U(N)$ 
if and only if 
$\hat{\chi}_1(a_1)=0$  
for $a_1\notin U(R_1)$, and
$\hat{\chi}_2(a_2)=0$ for  
$a_2\notin U(R_2)$,
we have the next result.  \\ 

\vspace{2mm}

\noindent
{\bf Corollary 2} \ {\it $\chi$ is primitive if and only if $\chi_1$ and $\chi_2$ are primitive. }  \\

\noindent
We can now construct the primitive multiplicative characters on $\Bbb{Z}/N$, 
\begin{eqnarray}
N=p_1^{c_1}...p_r^{c_r}.
\end{eqnarray}

\noindent
Choose  
primitive multiplicative characters $\chi_j$ in $\Bbb{Z}/p_j^{c_j}$,
$0\leq j \leq r$. 
If $c_j=1$, we can choose any non-trivial character on $\Bbb{Z}/p$ for $\chi_j$.
Otherwise $c_j>1$, and by Theorem 2 we can choose any
multiplicative character on $\Bbb{Z}/p_j^{c_j}$ whose DFT does not vanish at $1$.
By Theorem 3
\begin{eqnarray}
\chi(a)=\chi_1(a_1)...\chi_r(a_r), \ \ \ a_j\equiv a \!\!\! \mod p_j^{c_j}, \ \ 0\leq j\leq r,
\end{eqnarray}
is a primitive multiplicative character on $\Bbb{Z}/N$, and 
every primitive multiplicative character on $\Bbb{Z}/N$ has this form. 
Usually the definition of primitive multiplicative character is given in terms of induced characters [4].
A primitive multiplicative character is a multiplicative character not induced from a lower order character.
Our definition follows the program of using the CRT to bootstrap general results from results on factors. \\

\noindent

\noindent
Suppose $\chi_1$ and $\chi_2$ are 
primitive multiplicative characters on $\Bbb{Z}/R_1$ and on $\Bbb{Z}/R_2$,
and $\chi$ is the corresponding
primitive multiplicative character on $\Bbb{Z}/N$, $N=R_1R_2$, $(R_1, R_2)=1$.
Set
\begin{eqnarray}
f_1=\frac{e_1}{R_2} \ \ {\rm and} \ \
f_2=\frac{e_2}{R_1}. 
\end{eqnarray}
By primitivity 
\begin{eqnarray}
& \hat{\chi}_1(a_1f_1)
=G(\chi_1)\chi_1^{\ast}(a_1f_1)
=G(\chi_1)\chi_1^{\ast}(a_1)\chi_1^{\ast}(f_1), & \\
& \hat{\chi}_2(a_2f_2)
=G(\chi_2)\chi_2^{\ast}(a_2f_2)
=G(\chi_2)\chi_2^{\ast}(a_2)\chi_2^{\ast}(f_2), &
\end{eqnarray}
and
\begin{eqnarray}
\hat{\chi}(a_1e_1+a_2e_2)=G(\chi)\chi_1^{\ast}(a_1)\chi_2^{\ast}(a_2). \ \ a_1\in \Bbb{Z}/R_1, \ \ a_2\in \Bbb{Z}/R_2.
\end{eqnarray}
This leads to the next result. \\

\noindent
{\bf Theorem 4} \ {\it  
If $\chi_1$ and $\chi_2$ are 
primitive multiplicative characters on $\Bbb{Z}/R_1$ and on $\Bbb{Z}/R_2$,
and
$\chi$ is the corresponding primitive multiplicative character on $\Bbb{Z}/N$,
$N=R_1R_2$, $(R_1, R_2)=1$, then}
\begin{eqnarray}
\hat{\chi}(a_1e_1+a_2e_2)=G(\chi)\chi_1^{\ast}(a_1)\chi_2^{\ast}(a_2), \ \ a_1\in \Bbb{Z}/R_1, \ \ a_2\in \Bbb{Z}/R_2,
\end{eqnarray}
{\it where}
\begin{eqnarray}
G(\chi)=G(\chi_1)G(\chi_2)\chi_1^{\ast}(f_1)\chi_2^{\ast}(f_2).
\end{eqnarray}

\noindent
\section{DFT of multiplicative characters and the GTPF algorithm}

\noindent
In this section we derive the main result of Section 3, Theorem 4, in a slightly different way, using the GTPF algorithm. 
The GTPF algorithm was first described in [18] and [30]. It is an analoque of the Cooley-Tukey (CT) algorithm for computing the DFT.
While the CT algorithm addresses the case when $N=2^n$, the GTPF algorithm addresses the case when $N$ is an odd prime. 
The key components of the GTPF algorithm are the DFTs of the factors and permutation matrices. 
The details of the algebraic structure of the algorithm are given in [31]. \\

\noindent
Suppose $\chi_1$ and $\chi_2$ are primitive multiplicative characters on $\Bbb{Z}/R_1$ and $\Bbb{Z}/R_2$, and
$\chi$ is the corresponding primitive multiplicative characters on $\Bbb{Z}/N$. Then
\begin{eqnarray}
\hat{\chi}_1=G(\chi_1)\chi_1^{\ast}, \ \ \ 
\hat{\chi}_2=G(\chi_2)\chi_2^{\ast},  
\end{eqnarray} 
and
\begin{eqnarray}
\hat{\chi}=G(\chi)\chi^{\ast}.  
\end{eqnarray} 
We will now use the GTPF algorithm to relate the constants $G(\chi_1)$, $G(\chi_2)$ and $G(\chi)$. \\

\noindent
Every $a\in \Bbb{Z}/N$ can be uniquely written as
\begin{eqnarray}
a= a_2+a_1 R_2,  \ \ 0\leq a_1<R_1, \ \ 0\leq a_2<R_2.  
\end{eqnarray}
with $a_2$ the fastest running variable for consistency with the tensor product. 
%
Define the permutation $\pi$ of
 ${\Bbb{Z}}/N$ by 
\begin{eqnarray}
\pi(a)\equiv a_1e_1+a_2e_2 \!\!\! \mod N, 
\end{eqnarray} 
and set $Q_{\pi}$ equal to the corresponding permutation matrix. 
Then [31] 
\begin{eqnarray}
F(N)=Q_{\pi}^{-1} F_{\pi} Q_{\pi},
\end{eqnarray}
where
\begin{eqnarray}
F_{\pi}=
[w^{\pi(a)\pi(b)}]_{0\leq a,b<N}, \ \ w=e^{2\pi i\frac{1}{N}}.
\end{eqnarray}  

\noindent
We describe the structure of $F_{\pi}$ below. \\ 

\noindent
Applying $Q_{\pi}$ to $\chi(a)$ we have
\begin{eqnarray}
Q_{\pi}\chi(a)=\chi(\pi(a))=\chi(a_1e_1+a_2e_2)=\chi_1(a_1)\otimes \chi_2(a_2).
\end{eqnarray}
Set
\begin{eqnarray}
u_1=e^{2\pi i\frac{1}{R_1}}, \ \ \
u_2=e^{2\pi i\frac{1}{R_2}},
\end{eqnarray}
and
\begin{eqnarray}
f_1=\frac{e_1}{R_2}, \ \ \
f_2=\frac{e_2}{R_1}.
\end{eqnarray}
$f_1$ and $f_2$ are integers with $f_1$ relatively prime to $R_1$,
and $f_2$ relatively prime to $R_2$. Then
\begin{eqnarray}
f_1 t, \ \ \ 0\leq t <R_1,
\end{eqnarray}
is a permutation of
 ${\Bbb{Z}}/R_1$ and 
\begin{eqnarray}
f_2 t, \ \ \ 0\leq t <R_2
\end{eqnarray}
is a permutation of
 ${\Bbb{Z}}/R_2$. \\ 

\noindent
We can now describe $F_{\pi}$. 
From
\begin{eqnarray}
w^{\pi(a)\pi(b)}=u_1^{f_1a_1b_1}u_2^{f_2a_2b_2},
\end{eqnarray}
we have
\begin{eqnarray}
F_{\pi}=F_{R_1}\otimes F_{R_2}
\end{eqnarray}

\noindent
where
\begin{eqnarray}
F_{R_1}=
[u_1^{f_1a_1b_1}]_{0\leq a_1,b_1<R_1} 
\end{eqnarray}
and
\begin{eqnarray}
F_{R_2}=
[u_2^{f_2a_2b_2}]_{0\leq a_2,b_2<R_2}. 
\end{eqnarray}  
$F_{R_1}$ is found by permuting the rows of $F(R_1)$ by $f_1t, \ 0\leq t<R_1$,
with a similar statement for $F_{R_2}$. 
Then
\begin{eqnarray}
F_{\pi}Q_{\pi}\chi(a)= 
F_{R_1}(\chi_1(a_1))\otimes 
F_{R_2}(\chi_2(a_2))=G(\chi_1)G(\chi_2) \ 
\chi_1^{\ast}(f_1)\chi_2^{\ast}(f_2) 
\ \chi_1^{\ast}(a_1)\otimes \chi_2^{\ast}(a_2) 
\end{eqnarray}
and
\begin{eqnarray}
\hat{\chi}(a)=Q_{\pi}^{-1}F_{\pi}Q_{\pi}\chi(a)= 
G(\chi_1)G(\chi_2) \ 
\chi_1^{\ast}(f_1)\chi_2^{\ast}(f_2) 
\ \chi^{\ast}(a). 
\end{eqnarray}

\noindent
This leads to the next result. \\


\noindent
{\bf Theorem 5} \ {\it Suppose $\chi_1$ and $\chi_2$ are primitive multiplicative characters on $\Bbb{Z}/R_1$ and on $\Bbb{Z}/R_2$, and
\begin{eqnarray}
\chi(a)=\chi_1(a_1)\chi_2(a_2), \ \ \ a_1\equiv a \!\!\! \mod R_1 \ \ \ {\rm and} \ \ \ a_2\equiv a \!\!\! \mod R_2,
\end{eqnarray} 
the corresponding primitive multiplicative characters on $\Bbb{Z}/N$. 
Then}
\begin{eqnarray}
\hat{\chi}(a)=G(\chi_1)G(\chi_2)
\chi_1^{\ast}(f_1)\chi_2^{\ast}(f_2) 
\chi^{\ast}(a), \ \ a\in \Bbb{Z}/N. 
\end{eqnarray} 

\noindent
In the next two sections we use results from Sections 3 and 4 to derive the DFTs of the Jacobi symbol,
the modified Jacobi sequence and the Golomb sequence. \\

\section{Application of the DFT formula to the Jacobi symbol}

\noindent
Let $N=R_1R_2$, $(R_1,R_2)=1$,
and $\chi_{R_1}$ and $\chi_{R_2}$ be the Jacobi symbols on 
$\Bbb{Z}/R_1$
and
$\Bbb{Z}/R_2$, 
viewed as primitive multiplicative characters on
$\Bbb{Z}/R_1$
and
$\Bbb{Z}/R_2$ (Appendix A1), 
\begin{eqnarray}
\chi_{R_1}(a_1)=\left(\frac{a_1}{R_1}\right) \ \ \ {\rm and} \ \ \ 
\chi_{R_2}(a_2)=\left(\frac{a_2}{R_2}\right). 
\end{eqnarray} 
Then
\begin{eqnarray}
\left(\frac{R_2}{R_1}\right) 
\chi_{R_1}(f_1)=
\left(\frac{R_2}{R_1}\right) 
\left(\frac{e_1/R_2}{R_1}\right)= 
\left(\frac{e_1}{R_1}\right)=1 
\end{eqnarray} 
and
\begin{eqnarray}
\left(\frac{R_1}{R_2}\right) 
\chi_{R_2}(f_2)=
\left(\frac{R_1}{R_2}\right) 
\left(\frac{e_2/R_1}{R_2}\right)= 
\left(\frac{e_2}{R_2}\right)=1. 
\end{eqnarray} 

\vspace{2mm}

\noindent
By Theorem 4 
we have the next result. \\

\noindent
{\bf Theorem 6} \ {\it If $\chi_{R_1}$ and $\chi_{R_2}$ are the Jacobi symbols on $\Bbb{Z}/R_1$ and $\Bbb{Z}/R_2$,
and $\chi_N$, $N=R_1R_2$, $(R_1,R_2)=1$, is the corresponding Jacobi symbol on $\Bbb{Z}/N$, then}
\begin{eqnarray}
G(\chi_N)
=G(\chi_{R_1})G(\chi_{R_2})
\left(\frac{R_1}{R_2}\right)
\left(\frac{R_2}{R_1}\right).
\end{eqnarray} 

\vspace{4mm}

\noindent
Applying QRL in Appendix 1, Theorem 6 implies the following result. \\

\noindent
{\bf Corollary 3} \ {\it }
\begin{eqnarray}
G(\chi_N)
=G(\chi_{R_1})G(\chi_{R_2})
(-1)^{\frac{R_1-1}{2}\frac{R_2-1}{2}}.
\end{eqnarray} 

\vspace{2mm}

\noindent
The Legendre symbol 
\begin{eqnarray}
\chi_p(a)=
\left(\frac{a}{p}\right), 
\end{eqnarray} 
where $p$ is an odd prime, is a primitive multiplicative character on $\Bbb{Z}/p$.
We can identify $\chi_p$ with the unique real valued multiplicative character
\begin{eqnarray}
\chi_{t/2}, \ \ \ t=p-1
\end{eqnarray} 
in $\Bbb{Z}/p$. $\chi_p$ is an eigenvector of the $p$-point DFT,
\begin{eqnarray}
\hat{\chi}_p=c_p\chi_p, \ \ \ c_p=G(\chi_p).
\end{eqnarray} 

\noindent
By Theorem 4, the Jacobi symbol
\begin{eqnarray}
\chi_N(a)=
\left(\frac{a}{N}\right)= 
\left(\frac{a}{p_1...p_r}\right)= 
\left(\frac{a}{p_1}\right)... 
\left(\frac{a}{p_r}\right) 
\end{eqnarray} 
is a primitive multiplicative character on $\Bbb{Z}/N$,
\begin{eqnarray}
\hat{\chi}_N(a)=c_N\chi_N(a). 
\end{eqnarray} 
Applying Corollary 3 permits computing an explicit formula for $c_N$, which leads to the next result. \\

\noindent
{\bf Corollary 4} \ {\it
$\chi_N$ is an eigenvector of the $N$-point DFT, i.e.,
\begin{eqnarray}
\hat{\chi}_N(a)=c_N\chi_N(a), 
\end{eqnarray} 
where
\begin{eqnarray}
c_N=c_{p_1...p_r}=c_{p_1}...c_{p_r}
(-1)^{\frac{1}{4}\sum_{l=1}^{r-1}(p_1...p_{l}-1)(p_{l+1}-1)}, 
\end{eqnarray} 
and
\begin{eqnarray}
c_{p_i}=G(\chi_{p_i})=
\begin{cases}
\ \sqrt{p_i},     
&  p_i\equiv 1 \mod 4, \\ 
\ i\sqrt{p_i}, & p_i\equiv 3 \mod 4.
\end{cases} 
\end{eqnarray} 
}

\noindent
In particular, when $N=p_1p_2$,
\begin{eqnarray}
c_{p_1p_2}=c_{p_1}c_{p_2}
(-1)^{\frac{p_1-1}{2}\frac{p_2-1}{2}} 
\end{eqnarray} 
and, when $N=p_1p_2p_3$,
\begin{eqnarray}
c_{p_1p_2p_3}=c_{p_1}c_{p_2}c_{p_3}
(-1)^{\frac{p_1-1}{2}\frac{p_2-1}{2}} 
(-1)^{\frac{p_1p_2-1}{2}\frac{p_3-1}{2}}. 
\end{eqnarray} 

\noindent

\vspace{3mm}

\section{Application of the DFT formula to the modified Jacobi sequence and the Golomb sequence}

\noindent
Suppose $N=pq$, $p$ and $q$ are distinct odd primes, $p<q$.
Define the {modified} Legendre and Jacobi symbols, 
$x_p$ and $x_{pq}$, 
by\footnote{In sequence design literature, including in [22], $x_p(0)$ and $x_{pq}(0)$ are 
usually set to $1$. We choose here to follow the standard mathematical convention.}
\begin{eqnarray}
x_p(n):=
\begin{cases}
\ 0,     
&  n=0, \\ 
\ \rho\left( \left(\frac{n}{p}\right)\right), & {\rm otherwise}. 
\end{cases}
\end{eqnarray} 
and
\begin{eqnarray}
x_{pq}(n):=
\begin{cases}
\ 0,     
&  n=0, \\  
\ 1,     
&  n\equiv 0\mod p, \ n\not\equiv 0\mod q, \\ 
\ 0,     
&  n\equiv 0\mod q, \ n\not\equiv 0\mod p, \\ 
\ \rho\left(
\left(\frac{n}{p}\right)     
\left(\frac{n}{q}\right)
\right),     
&  (n,pq)=1,  
\end{cases}
\end{eqnarray} 
where $\rho(-1)=1$ and $\rho(1)=0$. 
The modified Legendre and Jacobi sequences are of interest in sequence design applications, as they
have large linear complexity 
and large merit factor [19], [22], [33].
They are also linked to the Golomb sequences [15],
a well-known unimodular two-valued ideal sequence. 
The DFT of the Golomb sequence, which we obtain by first computing the DFT of the modified Jacobi sequence,
 is the main result of this section. \\

\noindent
The modified Jacobi sequence in (88) can be re-written as
\begin{eqnarray}
x_{pq}(n)=
\frac{1}{2}\left(
\left(\frac{n}{pq}\right)
+1_{pq}-comb_p+comb_q+\delta_0\right),
\end{eqnarray} 
where $1_{pq}$, $comb_p$ and $\delta_0$ are the length-$pq$ vectors 
\begin{eqnarray}
1_{pq}=[1, 1, ... , 1],
\end{eqnarray} 
\begin{eqnarray}
comb_p=
\begin{cases}
\ 1,
&  
n\equiv 0\mod p, \\
\ 0, 
&  
{\rm otherwise}, 
\end{cases}
\end{eqnarray} 
and
\begin{eqnarray}
\delta_0=[1, 0, ... , 0].
\end{eqnarray} 
Since
\begin{eqnarray}
\hat{1}_{pq}=pq\delta_0,
\end{eqnarray} 
\begin{eqnarray}
\widehat{comb}_{p}=
\sum_{r=0}^{q-1}e^{2\pi i\frac{mr}{q}}=
\begin{cases}
\ q,
&  
m\equiv 0 \mod q, \\
\ 0, 
&  
{\rm otherwise}, 
\end{cases}
\end{eqnarray} 
\begin{eqnarray}
\hat{\delta}_0=1_{pq}, 
\end{eqnarray} 
and
\begin{eqnarray}
\hat{\left(\frac{n}{pq}\right)}=c_{pq}\left(\frac{m}{pq}\right)= 
\begin{cases}
\ 0,
&  
(m,pq)>1, \\ 
\ c_{pq}\left(\frac{m}{pq}\right),
&  
{\rm otherwise}, 
\end{cases}
\end{eqnarray} 
where
\begin{eqnarray}
c_{pq}=
\begin{cases}
\ i\sqrt{pq},
&  
q=p+2, \\ 
\ \sqrt{pq},
&  
{\rm otherwise}, 
\end{cases}
\end{eqnarray} 
which leads to the next result. \\

\noindent
{\bf Theorem 7}
\begin{eqnarray}
\hat{x}_{pq}(m) 
 =  \frac{1}{2}
\begin{cases}
\ pq-q+p+1,
&  
m=0, \\
\ -q+1,
&
m\equiv 0 \mod q, \ m\neq 0, \\
\ p+1,
&
m\equiv 0 \mod p, \ m\neq 0, \\
\ c_{pq}\left(\frac{n}{pq}\right)+1,
&  
{\rm otherwise}. 
\end{cases}
\end{eqnarray} 

\vspace{3mm}

\noindent
It follows that $|\hat{x}_{pq}(m)|$ is constant for all $m>0$ iff $q=p+2$. \\

\noindent
Define the Golomb sequence\footnote{
The Golomb sequence is
sometimes referred to as the Bj$\rm\ddot{o}$rck sequence [26], as it was discovered independently by Bj$\rm\ddot{o}$rck [5]. 
} 
as the following modification of the modified Jacobi sequence\footnote{
Golomb defines his sequence with respect to the Paley-Hadamard difference set, i.e.,
a cyclic difference set with 
the parameters $(n,k,\lambda)=(4t-1, 2t-1, t-1)$.
The construction of the cyclic Paley-Hadamard difference set is known for three cases: (1) $N=p$, (2) $=p(p+2)$, and (3) $N=2^m-1$.
For $q=p+2$ the modified Jacobi sequence is identical with the second case
(also known as the twin-prime sequence [29]), 
and for $N=p$ the modified Legendre sequence is identical with the first case.} 
[15],
\begin{eqnarray}
y_{pq}(n):=
\begin{cases}
\ 1,     
&  x_{pq}(n)=1, \\  
\ \alpha,     
&  x_{pq}(n)=0,  
\end{cases}
\end{eqnarray} 
where 
$\alpha=e^{i\Phi}$, $\Phi=cos^{-1}\left(-\frac{pq-1}{pq+1}\right)$ and $q=p+2$. \\

\noindent
The equation (99) can be re-written as
\begin{eqnarray}
y_{pq}=
(1-\alpha)x_{pq}+\alpha 1_{pq}. 
\end{eqnarray} 
Then
\begin{eqnarray}
\nonumber
\hat{y}_{pq}(m) & = &
\frac{1}{2}(1-\alpha)
\left(\hat{\left(\frac{n}{pq}\right)}-\widehat{comb}_p+\widehat{comb}_q+\hat{\delta}_0\right)
+\frac{1}{2}(1+\alpha)\hat{1}_{pq} \\
\end{eqnarray} 
This leads to the next result. \\

\noindent
{\bf Theorem 8}
\begin{eqnarray}
\hat{y}_{pq}(m) 
& = & \frac{1}{2}
\begin{cases}
\ (1+\alpha)pq+(1-\alpha)(-q+p+1),
&  
m=0, \\
\ (1-\alpha)(-q+1),
&
m\equiv 0 \mod q, \ m\neq 0, \\
\ (1-\alpha)(p+1),
&
m\equiv 0 \mod p, \ m\neq 0, \\
\ (1-\alpha)(c_{pq}\left(\frac{m}{pq}\right)+1),
&  
{\rm otherwise}.
\end{cases}
\end{eqnarray} 

\vspace{3mm}

\noindent
When $q=p+2$, $\hat{y}_{pq}$ 
has a constant magnitude 
and the Golomb sequence is ideal. \\

\noindent
When $N=p$, then
\begin{eqnarray}
x_{p}(n)=
\frac{1}{2}\left(
\left(\frac{n}{p}\right)
+1_{p}-\delta_0\right),
\end{eqnarray} 
\begin{eqnarray}
y_{p}=
(1-\alpha)x_{p}+\alpha 1_{p}, 
\end{eqnarray} 
\begin{eqnarray}
\hat{x}_{p}(m) 
 =  
\frac{1}{2}\left(\hat{\left(\frac{n}{p}\right)}+\hat{1}_{p}-\hat{\delta}_0\right) 
 =  \frac{1}{2}
\begin{cases}
\ p-1,
&  
m=0, \\
\ c_{p}\left(\frac{m}{p}\right)-1,
&  
{\rm otherwise}, 
\end{cases}
\end{eqnarray} 
and
\begin{eqnarray}
\hat{y}_{p}(m)=(1-\alpha)\hat{x}_p+\alpha \hat{1}_p 
 =  \frac{1}{2}
\begin{cases}
\ (1+\alpha)p+\alpha-1,
&  m=0, \\
\ (1-\alpha)(c_{p}\left(\frac{m}{p}\right)-1),
&  
{\rm otherwise}.
\end{cases}
\end{eqnarray} 

\noindent
When $p\equiv 3 \!\! \mod 4$, then $c_p=i\sqrt{p}$, $\hat{y}_p$ has a constant magnitude and the Golomb sequence is ideal.


\vspace{3mm}

\section{FZT of multiplicative characters}

\noindent
Suppose $N=LM$. Define the FZT of an arbitrary $N$-periodic sequence $x$ by
\begin{eqnarray}
X_L(j,k)=\sum_{r=0}^{L-1}x(k+rM)e^{2\pi i\frac{rj}{L}}, \ \ \ 0\leq j<L, \ \ 0\leq k<M.
\end{eqnarray}
FZT has several applications in mathematics [12], quantum
mechanics [34] and signal analysis [21]. In particular, it plays a major role in the analysis of time-frequency
representations, including ambiguity functions [3] and Weyl-Heisenberg expansions [1], 
and in polyphase sequence design [2], [6-11]. \\ 

\noindent
The FZT is a time-frequency representation that is closely related to the DFT [21].
For $L=1$ $X_L=x$. For $M=1$ $X_L=\hat{x}$. In general, $X_L$ can be viewed as a collection of DFTs of 
appropriate decimations of $x$ taken at various values of $k$.  
In particular, when $x$ is a discrete periodic chirp, the FZT of $x$ can be 
linked with the DFT of the decimation of an appropriate component of $x$ [10].
In effect, the condition for polyphase/semi-polyphase support of the DFT of a discrete periodic chirp
can be replaced by the condition of polyphase/semi-polyphase support of the FZT of a discrete periodic chirp, 
and the Fourier space design setting can be replaced with the Zak space design setting.
In this section we set the stage for future investigations of existence of similar relationships for binary sequences. \\

\noindent
Suppose $N=R_1R_2$, $(R_1, R_2)=1$, $\chi_1$ and $\chi_2$ are
multiplicative characters on $\Bbb{Z}/R_1$ and $\Bbb{Z}/R_2$, and $\chi$ is the corresponding
multiplicative character on $\Bbb{Z}/N$. 
In this section we describe the $R_1\times R_2$ FZT 
of $\chi$,
\begin{eqnarray}
X_{R_2}(j,k)=\sum_{r=0}^{R_2-1}\chi(k+rR_1)e^{2\pi i\frac{jr}{R_2}}.
\end{eqnarray} 
Suppose
\begin{eqnarray}
\chi(k+rR_1)=\chi_1(k)\chi_2(k+rR_1).
\end{eqnarray} 
We have
\begin{eqnarray}
X_{R_2}(j,k)=\chi_1(k)\sum_{r=0}^{R_2-1}\chi_2(k+rR_1)e^{2\pi i\frac{jr}{R_2}}, \ \ \ 0\leq j<R_2, \ \ 0\leq k<R_1.
\end{eqnarray} 
Consider
\begin{eqnarray}
\hat{\chi}_2(j)=\sum_{r=0}^{R_2-1}\chi_2(r)e^{2\pi i\frac{jr}{R_2}}.
\end{eqnarray} 
Multiplying by $\chi_2(R_1)$
\begin{eqnarray}
\chi_2(R_1)\hat{\chi}_2(j)=\sum_{r=0}^{R_2-1}\chi_2(rR_1)e^{2\pi i\frac{jr}{R_2}}.
\end{eqnarray} 
Changing variables by $rR_1=r'R_1+k$,
\begin{eqnarray}
\chi_2(R_1)\hat{\chi}_2(j)=e^{2\pi i \frac{jR_1^{-1}k}{R_2}}\sum_{r'=0}^{R_2-1}\chi_2(k+r'R_1)e^{2\pi i\frac{jr'}{R_2}}.
\end{eqnarray} 
Plugging this result into the formula for $X_{R_2}(j,k)$, we get
\begin{eqnarray}
X_{R_2}(j,k)=\chi_1(k)\chi_2(R_1)e^{-2\pi i \frac{jR_1^{-1}k}{R_2}}\hat{\chi}_2(j). 
\end{eqnarray} 
Suppose now that $\chi_2$ is primitive, then
\begin{eqnarray}
\hat{\chi}_2(j)=G(\chi_2)\chi_2^{\ast}(j), \ \ 0\leq j<R_2, 
\end{eqnarray} 
and we have the following result. \\

\noindent
{\bf Theorem 9} \ {\it If $\chi_1$ and $\chi_2$ are multiplicative characters on $\Bbb{Z}/R_1$ and $\Bbb{Z}/R_2$,
with $\chi_2$ primitive, then}
\begin{eqnarray}
X_{R_2}(j,k)=G(\chi_2)
\chi_2(R_1)
e^{-2\pi i \frac{R_1^{-1}jk}{R_2}}
\chi_2^{\ast}(j)
\chi_1(k), \ \ \
0\leq j<R_2, \ \ 0\leq k<R_1. 
\end{eqnarray} 

\vspace{4mm}

\noindent
Theorem 9 permits, in particular, an assessment of the DFT of the Jacobi symbol and related sequences. Set 
$\chi(k+rR_1)=
\left(\frac{k+rR_1}{R_1R_2}\right)$. We have the following result.  \\ 

\noindent
{\bf Corollary 5}
{\it
\begin{eqnarray}
X_{R_2}(j,k)=c_{R_2}
\left(\frac{R_1}{R_2}\right) 
e^{-2\pi i \frac{R_1^{-1}jk}{R_2}}
\left(\frac{j}{R_2}\right) 
\left(\frac{k}{R_1}\right), \ \ \ 
0\leq j<R_2, \ \ 0\leq k<R_1, 
\end{eqnarray}
where $c_{R_2}$ is computed via formula (82). \\ \\
}

\noindent
The FZT of the modified Jacobi sequence can be written as
\begin{eqnarray}
X_q\{x_{pq}\}(j,k)  = 
\frac{1}{2}
X_q
\left\{\left(\frac{n}{pq}\right)+1_{pq}-comb_p+comb_q+\delta_0\right\}.  
\end{eqnarray}
Since
\begin{eqnarray}
X_q\{1_{pq}\}(j,k)=
\sum_{r=0}^{q-1}e^{2\pi i\frac{rj}{q}}=
\begin{cases}
\ q,
&
j=0, \ 0\leq k<p, \\ 
\ 0,
&
else,
\end{cases}
\end{eqnarray} 

\begin{eqnarray}
X_q\{comb_{p}\}(j,k)=
\begin{cases}
\ \sum_{r=0}^{q-1}e^{2\pi i\frac{rj}{q}},
&
k=0, \\ 
\ 0,
&
else,
\end{cases}
=
\begin{cases}
\ q,
&
j=k=0, \\
\ 0,
&
else,
\end{cases}
\end{eqnarray} 

\begin{eqnarray}
X_q\{comb_{q}\}(j,k)=
\sum_{k+rp \ \equiv \ 0 \!\!\!\! \mod q}e^{2\pi i\frac{rj}{q}}=
e^{-2\pi i\frac{p^{-1}kj}{q}},
\ \ {\rm for \ all} \ \ j, k,
\end{eqnarray} 

\begin{eqnarray}
X_q\{\delta_0\}(j,k)=
\begin{cases}
\ 1,
&
k=0, \ \ 0\leq j<q, \\ 
\ 0,
&
else,
\end{cases}
\end{eqnarray} 
and 
\begin{eqnarray}
X_q\left\{\left(\frac{n}{pq}\right)\right\}(j,k)=
c_q 
\left(\frac{p}{q}\right)\left(\frac{k}{p}\right)\left(\frac{j}{q}\right)
e^{-2\pi i\frac{p^{-1}kj}{q}},
\ \ {\rm for \ all} \ \ j, k,
\end{eqnarray} 
we have the next result. \\

\noindent
{\bf Theorem 10}
{\it
\begin{eqnarray}
X_q\{x_{pq}\}(j,k) 
 =  
\begin{cases}
\ 2, &  k=0, \\ 
\ q+1, &  j=0, \ k\neq0, \\
\ A 
e^{-2\pi   i\frac{p^{-1}kj}{q}},  &  {\rm otherwise}, 
\end{cases}
\end{eqnarray} 
where
\begin{eqnarray}
A=1+c_q\left(\frac{p}{q}\right)\left(\frac{k}{p}\right)\left(\frac{j}{q}\right)
\end{eqnarray} 
and
\begin{eqnarray}
c_q=\sqrt{q}
\begin{cases}
\ 1, &  q\equiv 1 \!\!\! \mod 4, \\
\ i, &  q\equiv 3 \!\!\! \mod 4.
\end{cases}
\end{eqnarray} 
}

\noindent
Similarily, since the FZT of the Golomb sequence can be written as
\begin{eqnarray}
X_q\{y_{pq}\}  = 
(1-\alpha)X_q\{x_{pq}\}+\alpha X_q\{1_{pq}\},  
\end{eqnarray}
we have the next result. \\

\noindent
{\bf Theorem 11}
\begin{eqnarray}
X_q\{y_{pq}\}  =
\begin{cases} 
\ 1-\alpha+\alpha q,
&  
j=k=0, \\ 
\ \frac{(1-\alpha)(q+1)}{2}+\alpha q,
&  
j=0, \ k\neq 0, \\
\ 1-\alpha, 
&  
k=0, \ j\neq 0, \\
\ \frac{1-\alpha}{2}A 
e^{-2\pi i\frac{p^{-1}kj}{q}},  
&  
{\rm otherwise}.
\end{cases}
\end{eqnarray} \\

\vspace{-2mm}

\noindent
The array
$X_q\{x_{pq}\}(j>0,k>0)$ 
is unimodular when 
$q\equiv 3 \!\!\! \mod 4$, 
and semi-unimodular when
$q\equiv 1 \!\!\! \mod 4$.  
The support of 
$X_q\{y_{pq}\}$ 
is the same as the support of 
$X_q\{x_{pq}\}$.  \\ \\

\section{Further work}

\noindent
In this work we derived formulas for the $N$-point DFT and the $R_1\times R_2$ FZT of 
multiplicative characters on $\Bbb{Z}/N$, where $N=R_1R_2$ an odd integer
and $(R_1,R_2)=1$,
and of three special sequences: the Jacobi symbol, the modified Jacobi sequence and the Golomb sequence. \\ 

\noindent  
The importance of the results of this work 
stems from their potential utility in sequence design.
DFT and FZT enter 
sequence design at several stages in fundamental ways [2].
One of the main tasks in this field is identification of unimodular sequences with ideal periodic autocorrelation.
The ideal autocorrelation condition is equivalent to unimodularity of the DFT of a sequence.
Sequences satisfying this condition are called bi-unimodular or ideal [5], [26].
A corresponding condition exists for the FZT of a sequence [2], [11].
Other sequences of interest that are not ideal include Golay sequences [13] and zero autocorrelation zone sequences [17].
These sequences are often semi-unimodular, that is, their DFTs and FZTs have constant magnitude
on a subset of points and are zero otherwise.
The structure of the supports of the DFTs and the FZTs, 
have been previously used to design new families of chirp-like sequences [2], [6-11]. 
We anticipate the results for the DFT and FZT of multiplicative characters will be of equal utility. \\

\noindent
In particular,
these results 
avails convenient tools for investigation of two topics.
First, in the special case of Theorem 4, explored in Section 5, the DFT of the second order
primitive multiplicative characters was 
used to derive the DFT of 
the modified Jacobi sequence and the Golomb sequence. 
In further research it might be of interest to
consider the DFT of higher order primitive multiplicative characters, 
which might lead to identification of 
new sequences and description of their Fourier properties.  
Second, since, in general, there are many choices for the factors of $N$,
$R_1$ and $R_2$, there are many distinct, but equivalent factorizations of primitive
multiplicative characters on $\Bbb{Z}/N$, and hence many equivalent factorizations of the associated sequences.
The choice of a factorization or a subset of factorizations can be used in schemes
that require distributed processing or sharing of information between multiple users. \\

\noindent
Apart from sequence design applications, results of this work raise some broader questions:

\begin{enumerate}
\item
Primitive multiplicative characters are defined by the eigenvector property.
What is the relationship between eigenvectors and eigenvalues of different multiplicative characters?
\item
Modification of Jacobi sequences trades the eigenvector property for biunimodularity.
What is the mathematical relationship between the eigenvector property and biunimodularity?
Are there sequences other than ideal chirps, described in [6], that share these two properties?
\end{enumerate}

\noindent

\section*{Appendix A1 Synopsis of Legendre and Jacobi symbols}

\noindent
The Legendre symbol is defined by
\begin{eqnarray}
\left(\frac{a}{p}\right)
= 
\begin{cases}
\ 0,     
&  a=0, \\ 
\ 1,     
&  a \ {\rm is \ a \ quadratic \ residue} \!\!\! \mod p, \\  
\ -1,     
&  a \ {\rm is \ not \ a \ quadratic \ residue} \!\!\! \mod p. \\  
\end{cases}
\end{eqnarray}
The Legendre symbol is a primitive multiplicative character of order two on $\Bbb{Z}/p$ having
the following properties:
\begin{eqnarray}
& \left(\frac{ab}{p}\right)= 
\left(\frac{a}{p}\right) 
\left(\frac{b}{p}\right) & 
\end{eqnarray} 
\begin{eqnarray}
&{\rm if} \ a\equiv b \!\!\! \mod  p \ \ {\rm then} \ 
\left(\frac{a}{p}\right)= 
\left(\frac{b}{p}\right) & 
\end{eqnarray} 
\begin{eqnarray}
&\left(\frac{a^2}{p}\right)=1 & 
\end{eqnarray} 
In addition,
\begin{eqnarray}
a^{\frac{p-1}{2}}\equiv\left(\frac{a}{p}\right) \!\!\! \mod p. 
\end{eqnarray} 
If $p$ and $q$ are odd primes, the Legendre symbols satisfy the QRL,
\begin{eqnarray}
\left(\frac{p}{q}\right) 
\left(\frac{q}{p}\right)=(-1)^{\frac{p-1}{2}\frac{q-1}{2}}. \ \ \ \ 
\end{eqnarray} 
The Jacobi symbol is defined by
\begin{eqnarray}
\left(\frac{a}{N}\right)= 
\left(\frac{a}{p_1...p_r}\right)= 
\left(\frac{a}{p_1}\right) ...
\left(\frac{a}{p_r}\right). 
\end{eqnarray} 
Since the Legendre symbol has order two, we can assume that the primes $p_1$, ..., $p_r$ are distinct.
By Theorem 4 the Jacobi symbol is a primitive multiplicative character of order two on $\Bbb{Z}/N$ having the following properties:
\begin{eqnarray}
 \left(\frac{ab}{N}\right)= 
\left(\frac{a}{N}\right) 
\left(\frac{b}{N}\right) 
\end{eqnarray} 
\begin{eqnarray}
 {\rm if} \ \ a\equiv b \!\!\! \mod  N \ \ {\rm then} \ \ 
\left(\frac{a}{N}\right)= 
\left(\frac{b}{N}\right)  
\end{eqnarray} 
In addition,
\begin{eqnarray}
\left(\frac{a}{N_1N_2}\right)= 
\left(\frac{a}{N_1}\right) 
\left(\frac{a}{N_2}\right). 
\end{eqnarray} 
If $a$ and $b$ are odd positive integers, the Jacobi symbols satisfy the QRL,
\begin{eqnarray}
\left(\frac{a}{b}\right) 
\left(\frac{b}{a}\right)=(-1)^{\frac{a-1}{2}\frac{b-1}{2}}. 
\end{eqnarray}  


\vspace{3mm}


\section*{Appendix A2 Proof of the QRL}

\noindent
The QRL and the quadratic Gauss sum computation for the Jacobi symbol,
first derived by Gauss [20], are fundamental results in classical number theory.
These Gauss sums are critical to the determination of the class number of quadratic number fields,
as well as many other problems. Gauss also observed the connection between these sums and the QRL.
Since Gauss there have been many derivations of these results. \\ 

\noindent
In this appendix we describe another approach towards relating these results based on combining two independent formulas.
One, given in Theorem 6, uses the CRT to relate quadratic Gauss sums to the Jacobi symbol.
The second, derived in [13], 
uses the algebra of null-theta functions in the 3 dimensional real Heisenberg group
to relate the quadratic Gauss sum to the trace of the DFT matrix. \\

\noindent 

\noindent
We start by retracing the principal steps of the second approach. \\

\noindent
A bases for $L(\Bbb{Z}/p)$, the complex valued functions on $\Bbb{Z}/p$,
can be chosen, so that in this bases
\begin{eqnarray}
\frac{1}{\sqrt{p}}
F(p) 
\end{eqnarray} 
is a block diagonal matrix consisting of $2\times 2$ blocks having
zero trace and the single eigenvalue
\begin{eqnarray}
\frac{1}{\sqrt{p}}
G(\chi_p)
\end{eqnarray}
corresponding to the eigenvector $\chi_p$. Then 
\begin{eqnarray}
\frac{1}{\sqrt{p}}
G(\chi_p)=
Tr\left(\frac{1}{\sqrt{p}}F(p)\right)
=\frac{1}{\sqrt{p}}
\sum_{x=0}^{p-1}e^{2\pi i\frac{x^2}{p}}. 
\end{eqnarray} 
Details can be found in [31]. \\

\noindent
In the same way
\begin{eqnarray}
\frac{1}{\sqrt{N}}
G(\chi_N)=
Tr\left(\frac{1}{\sqrt{N}}F(N)\right)=
\frac{1}{\sqrt{N}}
\sum_{x=0}^{N-1}e^{2\pi i\frac{x^2}{N}}. 
\end{eqnarray}

\noindent 

\vspace{2mm}

\noindent
In [13] we derived a formula for
$Tr(\frac{1}{\sqrt{N}}F(N))$, independent of the QRL, 
based on the algebra of nil-theta functions on the 3 dimensional Heisenberg group,
\begin{eqnarray}
\frac{1}{\sqrt{N}}
G(\chi_N)=
Tr\left(\frac{1}{\sqrt{N}}F(N)\right)=
\begin{cases}
\ 1+i,     
&  N\equiv 0 \!\!\! \mod 4, \\  
\ 1,     
&  N\equiv 1 \!\!\! \mod 4, \\  
\ 0,     
&  N\equiv 2 \!\!\! \mod 4, \\  
\ i,     
&  N\equiv 3 \!\!\! \mod 4.  
\end{cases}
\end{eqnarray} 
In particular, when $N$ is odd,
\begin{eqnarray}
\frac{1}{\sqrt{N}}
G(\chi_N)=
i^{\frac{N-1}{2}}.
\end{eqnarray} 
Since, by (145),
\begin{eqnarray}
\frac{G(\chi_N)}{G(\chi_{R_1})G(\chi_{R_2})}=
i^{\frac{N-R_1-R_2+1}{2}},
\end{eqnarray} 
we have
\begin{eqnarray}
G(\chi_N)=G(\chi_{R_1})G(\chi_{R_2})(-1)^{\frac{R_1-1}{2}\frac{R_2-1}{2}}. 
\end{eqnarray}
Combining (147) with Theorem 6, 
\begin{eqnarray}
G(\chi_N)=G(\chi_{R_1})G(\chi_{R_2})
\left(\frac{R_1}{R_2}\right)
\left(\frac{R_2}{R_1}\right),
\end{eqnarray}
we have the QRL for the Jacobi symbol,
\begin{eqnarray}
\left(\frac{R_1}{R_2}\right)
\left(\frac{R_2}{R_1}\right)
=(-1)^{\frac{R_1-1}{2}\frac{R_2-1}{2}}. 
\end{eqnarray}

\end{document}